\documentclass[a4paper,12pt,draft]{amsart}

\usepackage[cp1252]{inputenc}
\usepackage{stmaryrd}
\usepackage{amsthm}
\usepackage{amscd}
\usepackage[all]{xy}

\begin{document}

\newtheorem{THM}{{\!}} [section]
\newtheorem{THMX}{{\!}} 
\renewcommand{\theTHMX}{}
\newtheorem{thm}{Theorem} [section]
\newtheorem{cor}[thm]{Corollary}
\newtheorem{lem}[thm]{Lemma}
\newtheorem{prop}[thm]{Proposition}
\newtheorem{defi}[thm]{Definition}
\newtheorem{axiom}[thm]{Axiom}
\newtheorem{rem}[thm]{Remark}
\newtheorem{souv}[thm]{Reminder}
\newtheorem{ex}[thm]{Example}
\newtheorem{exo}[thm]{Exercice}
\newtheorem{dem}[thm]{Proof}

\newcommand{\Z}{\mathbb{Z}}
\newcommand{\N}{\mathbb{N}}
\newcommand{\Q}{\mathbb{Q}}
\newcommand{\R}{\mathbb{R}}


\title[Subanalytic sets of low arity]{{\it o}-minimal structures:\\ low arity versus generation} 

\author{Serge Randriambololona}
\address{Université de Savoie, Le Bourget-du-Lac, France}
\email{Serge.Randriambololona@ens-lyon.org}

\date{March 31, 2005}


\begin{abstract}

We show that an analogue of the Hilbert's Thirteenth Problem fails in the real subanalytic setting.
Namely we prove that, for any integer $n$, the $o$-minimal structure generated by restricted analytic functions 
in $n$ variables is strictly smaller than the structure of all global subanalytic sets, whereas these two structures 
define the same subsets in $\R ^{n+1}$.  

\end{abstract}

\maketitle





%
%

\section{Introduction}

The aim of this paper is to prove that, for any fixed $n\in \N$, 
the $o$-minimal structure generated by the family of all  global subanalytic subsets of $\R ^n$
is strictly smaller than the structure of all global subanalytic sets: 
some subanalytic subsets of $\R ^{n+1}$ are ``transcendental'' 
over the family of all subanalytic subsets of $\R ^n$.

The main motivation for this work was to prove that the statement \begin{quote}
``{\it Given an $o$-minimal structure 
$\mathcal{S}$ over $X$, there is an integer $n$ such that 
$\mathcal{S}$ and $\mathrm{str}\, (\mathcal{S}^{(n)})$ - its reduct generated by $\mathcal{S}$-definable 
subsets of $X^n$ - define the same subsets of $X^N$, for all $N$}''
\end{quote} 
is false. We now know it fails for $\mathcal{S}$ being 
the structure of global subanalytic sets.

This result can be seen as a negative answer to a generalized real analytic version of the second part of 
Hilbert's Thirteenth Problem: subanalytic functions don't have 
the superposition property (see \cite{Kol} for the positive answer in the continuous setting).


In section \ref{definition}, we give the following definitions: $o$-minimal structure, generated structure, 
subanalytic sets and sub-$n$-analytic sets; only the last one is original. 
We then recall some well known properties.
 
In section \ref{Dan} is proven that restricted analytic functions in $n$ variables and 
subanalytic subsets of $\R ^{n+1}$ have the same definability power.
This elegant proof is due to Daniel J. Miler and is based on Hironaka's Uniformization Theorem 
for subanalytic sets.
  
In section \ref{subreg}-\ref{construction}, we use Gabrielov's 
``Explicit Fibre Cutting Lemma'', a diagonal argument on formal series and metric control on truncation of translated 
power series, to prove that there is a restricted analytic 
function $f:[-1,1]^{n+1} \to \R$ whose graph can't be defined by mean of
restricted analytic functions in $n$ variables.

\section{Definitions} 
\label{definition}

\begin{defi}

We call $\mathcal{S}=(\mathcal{S}^{(n)})_{n\in \N}$ a structure over $(\R;+,\cdot\, )$ if it has the following properties
\begin{itemize}
\item[(S1)] \label{Boole} $\mathcal{S}^{(n)}$ is a boolean subalgebra of $\mathcal{P}(\R ^n)$ for each $n\in \N$,
\item[(S2)] \label{SemAlg} if $n$ is an integer and $A$ is a semialgebraic subset of $\R ^n$ then $A \in \mathcal{S}^{(n)}$,
\item[(S3)] \label{Cart} if $A\in \mathcal{S}^{(n)}$, then $\R \times A\in \mathcal{S}^{(n+1)}$
\item[(S4)] \label{Proj} if $A \in \mathcal{S}^{(n+1)}$ and $\pi: \R^{n+1} \to \R ^n$ is the cartesian projection 
$\pi (x_1, \ldots , x_{n+1}) = (x_1, \ldots , x_n)$ then $\pi (A) \in \mathcal{S}^{(n)}$

\end{itemize}
If it furthermore has the property:
\begin{itemize}\item[(S5)] \label{omin} every element of ${\mathcal S}^{(1)}$ is a finite union of singletons and open intervals,
\end{itemize}
it is said to be an $o$-minimal structure over $(\R;+,\cdot\,)$.
\end{defi}

In words, a structure over $(\R;+,\cdot)$ is a collection of real sets, containing the family of all semialgebraic sets 
and stable under natural set theoretical operations: union, intersection, complementation, 
cartesian projection and cartesian product. The structure is $o$-minimal if the elements of $\mathcal{S}^{(1)}$ are the 
simplest possible: finite union of intervals and points.

\smallskip 

Elements of $\bigcup _n \mathcal{S}^{(n)}$ are called $\mathcal{S}$-definable sets; given a $\mathcal{S}$-definable set $A$, 
we call the integer $n$ such that $A \in \mathcal{S}^{(n)}$ the \emph{arity} of $A$.

A function $f$ from some $A\subseteq \R ^n$ to $\R ^m$ is said to be $\mathcal{S}$-definable if its graph is a 
$\mathcal{S}$-definable set. 

For an introduction to the geometry in $o$-minimal structure, see for instance \cite{Co1} or 
\cite{vdD}. 

\smallskip

Let us now define the notion of \emph{generated structure}.

If $\mathcal{U}=(\mathcal{U}^{(n)})_{n\in \N}$ and $\mathcal{V}=(\mathcal{V}^{(n)})_{n\in \N}$ are such 
that $\mathcal{U}^{(n)}\subseteq \mathcal{P} (\R ^n)$ and $\mathcal{V}^{(n)}\subseteq \mathcal{P} (\R ^n)$ 
, we will note $\mathcal{U} \sqsubseteq \mathcal{V}$  
the property `` $\mathcal{U} ^{(n)} \subseteq \mathcal{V} ^{(n)}$ for all $n \in \N$ ''.

If $\mathcal{A}=(\mathcal{A} ^{(n)})_{n\in \N}$ is such that $\mathcal{A} ^{(n)} \subseteq \mathcal{P}(\R ^n)$, there exists 
a smallest element - for the partial order $\sqsubseteq$ on $\prod_{n\in \N} \mathcal{P} (\mathcal{P} (\R^n))$ - 
among the $\mathcal{S}=(\mathcal{S}^{(n)})_{n\in \N}$ forming a structure over $(\R;+,\cdot\, )$   
and satisfying $\mathcal{A}\sqsubseteq \mathcal{S}$. We will 
 note this structure $\mathrm{str}\, (\mathcal{A})$, and call it the 
\emph{structure generated by $\mathcal{A}$}. 

\begin{rem}
Let $n_0$ be an integer and $\mathcal{F}^{(n_0)}$ a subset of $\mathcal{P}(\R ^{n_0})$; 
we will, when no confusion is possible, identify $\mathcal{F}^{(n_0)}$ and the family
$$\mathcal{G}=(\mathcal{G}^{(n)})_{n\in \N} \in \prod _{n\in \N} \mathcal{P}(\mathcal{P}(\R ^n)),$$ where 
$\mathcal{G}^{(n)}=\emptyset$ if 
$n\neq n_0$ and $\mathcal{G}^{(n_0)} = \mathcal{F}^{(n_0)}$. 

In such a case $\mathrm{str}\, (\mathcal{F}^{(n_0)})$ stands 
for $\mathrm{str}\, (\mathcal{G})$.
\end{rem}

Given an $n\in \N$, we let ${\mathcal B}(n)$ be the algebra of all functions $f:[-1,1]^n \to \R$
such that $f$ admits an analytical continuation in a neighbourhood of $[-1,1]^n$. We call such a function $f$ a 
\emph{restricted analytic function} (in $n$ variables).

Let $\mathcal{E}=(\mathcal{E}^{(n)})_{n \in \N ^*}$ be the element of $\prod _{n\in \N ^*} \mathcal{P}(\mathcal{P}(\R ^n))$
defined by $$\mathcal{E}^{(n+1)} :=\{ \mathrm{graph}(f) \, , \, f\in \mathcal{B}(n)\}.$$

With the previous notation, we denote by $\R _ {\mathrm{an}} $ the structure $\mathrm{str}\,( \mathcal{E} )$.

\begin{thm}[Gabrielov]
$\R _ {\mathrm{an}}$ is an $o$-minimal structure.
\end{thm}

An element $A$ in $\R_{\mathrm{an}}$ is called a \emph{global subanalytic set}. 

\begin{defi}
Given an integer $n$ we let $$\R_{\mathrm{an}(n)}:=\mathrm{str}\,(\mathcal{E}^{(n+1)});$$ 
$\R_{\mathrm{an}(n)}$-definable sets are called \emph{global sub-$n$-analytic sets}. 
\end{defi}
In words, $\R_{\mathrm{an}(n)}$ is the structure generated by the graphs of all restricted analytic functions \emph{in 
at most $n$ variables} 
(whereas there is no bound on the number of variables for the restricted analytic functions used to generate $\R_{\mathrm{an}}$).

For instance, $$\{ (x_1,x_2,x_3)\in [-1,1]^3 ; \cos \frac{x_1+ x_2}{2} + \sin \frac{x_3 - \cos x_2}{2}>0 \}$$ is a 
$\R_{\mathrm{an}(1)}$-definable subset of $\R ^3$.

\begin{prop}
\label{toto}
$\R_{\mathrm{an}(n)}$ is model complete (as a $\mathcal{B}(n)$-structure).
\end{prop}

Let $p$ be an integer; we will denote by $A_{p}({\mathcal B}(n))$ the subalgebra  of 
${\mathcal B}(p)$ generated by all the functions
$$ (x_1, \ldots , x_p) \mapsto f(x_{\sigma(1)} , \ldots , x_{\sigma(n)}),$$ 
as $\sigma$ ranges in $\{1,\ldots , p \}^{\{1,\ldots , n \}}$ (the set of functions from ${\{1,\ldots , n \}}$ to 
$\{1,\ldots , p \}$)
 and $f$ ranges in ${\mathcal B}(n)$ (the set of restricted analytic functions in $n$ variables).

Once we have noted that, for every $p\in \N$, the algebra  $A_{p}({\mathcal B}(n))$ is stable under the action of 
partial derivation operators, Proposition \ref{toto} easily follows from Gabrielov's ``Explicit Model Completeness'' 
(\cite{Gab}, Theorem 1. and Corollary).

We will use a more precise version of this result in sections \ref{subreg} and \ref{enumeration} to show how $\R_{\mathrm{an}(n)}$-definable 
functions are controlled by restricted analytic functions in at most $n$ variables.


\section{Sub-n-analytic sets}
\label{Dan}
\begin{prop}
\label{Miller}
$\R_{\mathrm{an}(n)}$ is the structure generated by global subanalytic sets of arity $n+1$.
\end{prop}

(Proof due to { Daniel J. Miller}.)

The inclusion $\R_{\mathrm{an}(n)}\sqsubseteq \mathrm{str}\,(\, \R_{\mathrm{an}}^{(n+1)})$ is easy.

Let us prove the other inclusion by induction on $n$; the case $n=0$ is clear.
Let denote by $K$ the set $[-1,1]^n$. By the cell decomposition theorem (\cite{vdD}, theorem 2.11.), 
it's enough to prove that, given a $\R_{\mathrm{an}}$-definable 
function $f:C\to \R$ for $C$ being a $\R_{\mathrm{an}}$-cell either included in or disjoint with $K$, then $f$
is $\R_{\mathrm{an}(n)}$-definable.  

Note that the mapping $i:(x_1, \ldots , x_n)\mapsto (1/x_1, \ldots , 1/x_n)$ 
is $\R_{\mathrm{an}(n)}$-definable and sends $\R \setminus K$ in $K$; we thus can suppose that $A\subseteq K$.

Up to a finer cell decomposition, we can furthermore suppose that $|f(\overline{x})|-1$ has constant sign on $C$ and, 
$y\mapsto 1/y$ 
being $\R_{\mathrm{an}(n)}$-definable, we can assume that $|f(\overline{x})|\leq 1$ for all $\overline{x}\in C$.

Let $G$ be the closure of the graph of $f$; $G$ is a compact subanalytic set of dimension $d \leq n$.

Hironaka's uniformization theorem (\cite{B-M}, theorem 0.1.) gives a $d$-dimensional analytic manifold $Y$ and a 
surjective analytic 
proper mapping $\psi:Y \to G$. 

$G$ being compact and $\psi$ being surjective and proper, $Y$ is compact; we then easily get a finite family $\{ \phi _i: [-1,1]^d \to Y \}_{i=1,\ldots , s }$ 
of restricted analytic functions such that the union of their images is covering $Y$. 

Hence $G = \bigcup _{i=1} ^s \psi \circ \phi _i ([-1,1]^d)$ is a $\R_{\mathrm{an}(d)}$-definable set and thus a 
$\R_{\mathrm{an}(n)}$-definable set.

By induction hypothesis, $C$ is an $\R_{\mathrm{an}(n-1)}$-definable set and thus a 
$\R_{\mathrm{an}(n)}$-definable set. 
The function $f$  is $\R_{\mathrm{an}(n)}$-definable, for its graph, $G \cap (C \times \R)$, is.


\section{n-regularity}

In the following sections, we prove that there are some $\R_{\mathrm{an}}$-definable analytic functions in $n+1$ variables 
which are not  $\R_{\mathrm{an}(n)}$-definable. We first show how each $\R_{\mathrm{an}(n)}$-definable function 
is ``controlled'', through the notion of 
\emph{ $n$-regularity}, by the restricted analytic functions in $n$ variables used to define it.

\label{subreg}
 Let $n$ and $p$ be two integers; in the proof of the Proposition \ref{toto}, we have defined 
 the algebra $A_{p}({\mathcal B}(n))$.

By definition, each $g\in A_{p}({\mathcal B}(n))$ can be written in the form 
$$g(x_1, \ldots, x_p) = Q\big( h_1 (x_{\sigma _1 (1)} , \ldots ,x_{\sigma _1 (n)}), \ldots ,  
h_q (x_{\sigma _q (1)} , \ldots ,x_{\sigma _q (n)})\big)$$
where $q$ is an integer, $Q$ is a polynomial in $q$ variables with integer coefficients, the $h_i$'s are restricted analytic functions in $n$ variables 
and the $\sigma _i$'s are mappings from $\{1, \ldots, n\}$ to $\{1, \ldots ,p \}$.

We will call an element of $A_{p}({\mathcal B}(n))$ a \emph{restricted analytic function in $p$ variables which essentially depends on at 
most $n$ variables}.

In some sense, the graph of a $\R_{\mathrm{an}(n)}$-definable function looks almost everywhere like an analytic variety defined as a zero-set of 
restricted analytic functions depending on at most $n$ variables.

Let's make this statement more precise: we first recall a special case of Gabrielov's ``Explicit Fibre Cutting Lemma'' 
(see \cite{Gab}, Lemma 3. and Theorem 1.):

\begin{thm}[Gabrielov]
\label{Gabrielov}
Given a $d$-dimensional sub-$n$-analytic set $Y\subseteq \R ^m$, there is a $p\in \N$,  
a finite family $\{X _{\nu} \}$ of sub-$n$-analytic subsets of $\R ^{m+p}$ and 
a sub-$n$-analytic subset $V$ of $\R ^{m+p}$ such that, if 
$\pi:\R ^m\times \R^p \to \R ^m$ is given by $\pi (x_1, \ldots , x_{m+p})=(x_1, \ldots , x_m)$, one has 
\begin{enumerate}
\item $Y= \pi (V) \cup \bigcup \pi (X_{\nu})$ ;
\item $\dim \pi (V) <d$ ;
\item for each $\nu$, $\dim X_{\nu} =d$ and $\pi _{|X_{\nu}} : X_{\nu} \to Y$ has rank $d$ at every point of $X_{\nu}$;
\item \label{varietelocale} for each $\overline{s}\in X_{\nu}$, $\{\overline{x}-\overline{s}\, ; \overline{x}\in X_{\nu} \}$ 
is near $\overline{0}$ the zero-set of $m+p-d$ elements 
$f_i :\R ^{m+p}\to \R$ of $A_{m+p}({\mathcal B}(n))$, $(df_i)_i$ having rank $m+p-d$ at $\overline{0}$;   
\item $X_{\lambda} \cap X_{\mu} = \emptyset$ for $\lambda \neq \mu$.
\end{enumerate}
\end{thm}

This theorem leads us to the following definition:

\begin{defi}
\label{subregular}
Let $f$ be a function, from a neighbourhood $U$ of $\overline{0}$ in $\R^{n+1}$, to $\R$.

$f$ is said to be $n$-regular at $\overline{0}$ if there exist 
\begin{itemize}
\item an integer $p$,
\item a $(p+1)$-tuple $(g_1, \cdots , g_{p+1})$ of elements of $A_{n+p+2}({\mathcal B}_n)$,
\item a neighbourhood $V\subseteq U$ of $\overline{0}\in \R^{n+1}$ and
\item  for each $\overline{x}\in V$, there is a point
$(y_1(\overline{x}) , \ldots, y_p(\overline{x}))$ in $\R ^p$, 
\end{itemize}

such that

\begin{itemize} 
\item $g_i(\overline{x}, y_1(\overline{x}), \ldots , y_p(\overline{x}),f (\overline{x}) ) =0$, for all $i$, and
\item the rank of $$\big( \frac{\partial g_i}{\partial z_j}\big) _{\substack{1\leq i \leq p+1 \\ n+2\leq j \leq n+p+2}}$$ 
is full at the point $(\overline{x}, y_1(\overline{x}) , \ldots , y_p(\overline{x}) ,f (\overline{x}))$. 
\end{itemize}

A function $f$ from a neighbourhood $U$ of $\overline{a} \in \R ^{n+1}$ is said to be $n$-subregular at $\overline{a}$ if 
$\overline{x} \mapsto f(\overline{a} + \overline{x})$ is  $n$-regular at $\overline{0}$.

\end{defi}

In words, $f$ is $n$-regular at $\overline{0}$ if, as in Theorem \ref{Gabrielov}, the germ of its graph is the germ of the
projection $\pi (X)$ of an analytic manifold $X$ given as the zero-set of some functions depending essentially on at most $n$ variables, 
and $\pi_{ |X}$ is locally a diffeomorphism.

\begin{prop}
Given a $\R_{\mathrm{an}(n)}$-definable function $f:[-1,1]^{n+1} \to \R$, there is a point $\overline{a} \in ]-1,1[^n$ such that 
$f$ is $n$-regular at $\overline{a}$.
\end{prop}

This proposition follows from an easy dimensional argument and Theorem \ref{Gabrielov}.

\section{Diagonalization}

\label{enumeration}
In the sequel we will build a function $h:[-1,1]^{n+1} \to \R$ such that 
\begin{itemize}
\item there is no $\overline{a} \in ]-1,1[^{n+1}$ at which $h$ is $n$-regular (and thus $h$ can't be  
$\R_{\mathrm{an}(n)}$-definable)
\item but $h$ is a restriction to $[-1,1]^{n+1}$ of some analytic function from $\R ^{n+1}$ to $\R$ (and subsequently is 
$\R_{\mathrm{an}}$-definable).
\end{itemize}


We will now ``enumerate'' the germs (above $\overline{0} \in \R ^{n+1}$) of $n$-regular (at $\overline{0}$) 
functions $f: \R ^{n+1} \to\R$.

We first have to choose a value $y$ for $f(0,\ldots, 0)$.

By definition of $n$-regularity, it is enough to look, as $p$ ranges in $\N$, at all the $(p+1)$-tuples $(g_1,\ldots, g_{p+1})$ 
of elements in $\mathcal{A} _{n+p+2} (\mathcal{B}_n)$ such that 
$g_i(0,\ldots,0,y) =0$ 
and the rank of $$\big( \frac{\partial g_i}{\partial z_j}\big) _{\substack{1\leq i \leq p+1 \\ n+2\leq j \leq n+p+2}}$$ 
is full at points $(0, \ldots ,0,y)$.

Let's fix such a $p\in \N$.

But by definition, each $g\in \mathcal{A} _{n+p+2} (\mathcal{B}_n)$ is of the following form:
\begin{itemize}
\item there is a $q\in \N$ and a $Q \in \Z [T_1, \ldots , T_q ]$,
\item there are some $h_1$, \ldots , $h_q$ in $\mathcal{B}(n)$,
\item for each $i\in \{1, \ldots ,q \}$, there is a mapping $\sigma _i$ from $\{1,\ldots ,n \}$ to $\{1, \ldots ,n+p+2 \}$
\end{itemize}
such that $$g(x_1, \ldots , x_{n+p+2})=Q(h_1 (x_{\sigma _1 (1)} , \ldots , x_{\sigma _1 (n)}), \ldots , 
h_q( x_{\sigma _q (1)} , \ldots , x_{\sigma _q (n)})).$$

So let's fix a $q\in \N$ a $(p+1)$-tuple of elements in $\Z [T_1, \ldots , T_q ]$ and for each $1\leq j \leq q$ and $1\leq i \leq p+1$, 
fix a mapping $\sigma _j ^i$ 
from $\{1,\ldots ,n \}$ to $\{1, \ldots ,n+p+2 \}$.
 
The only parameters left free 
are now 
\begin{itemize}
\item the value of $y$ of $f(0,\ldots, 0)$, 
\item the $(p+1)q$-tuple of restricted analytic functions $h$ in $n$ variables.
\end{itemize}

All those germs are thus built by choosing a set of ``assembly instructions'' (the integers $p$ and $q$, polynomials $Q$ and mappings $\sigma$) 
and then by assembling ``pieces'' (the restricted analytic functions $h$ in $n$ variables) that fit this set of instructions. 

Let
$$\left. \begin{array}{cccclll} 	
s		&\mapsto	
&	&\big((p(s)	&,q(s)),	&(Q_k (s))_{1 \leq k \leq p(s)+1}, 
			& (\sigma _j ^i (s) )_{\substack{1\leq j \leq q(s) \\ 1 \leq k \leq p(s)+1}}\big)
 \end{array} \right.$$ 
be a surjective mapping from $\N$ to $$\coprod _{\substack{(p,q)\in \N ^2}} \{ (p,q) \} \times {(\Z [T_1, \ldots , T_q])} ^{p+1}
\times{({(\{ 1,\ldots, n+p+2 \}^{\{ 1,\ldots, n\} })}^{q})}^{p+1}.$$

Fix an $s\in \N$ (and thus some integers $p(s)$, $q(s)$, some polynomials $(Q_k (s))_{1 \leq k \leq p(s)+1}$ and some mappings 
$(\sigma _j ^k (s) )_{\substack{1\leq j \leq q(s) \\ 1 \leq k \leq p(s)+1}}$). 

Then let $M_s$ be the subset of 
$$\R \times {(\R \{ X_1, \ldots , X_n \} ^{q(s)})}^{p(s)+1}$$ which elements 
$$\bigg(y,{\big( {(g^k _j)}_{1\leq j \leq q(s)}\big)}_{1\leq k \leq p(s) +1}\bigg)$$
satisfy the conditions in Definition \ref{subregular}:
\begin{enumerate}
\item $h_i(0,\ldots, 0,y)=0 , \forall i\in \{1 , \ldots , p(s)+1 \}$
\item the rank of $\displaystyle (\frac{\partial h_i}{\partial x_j})_{\substack{1\leq i \leq p+1 \\ n+1 \leq j \leq n+p+2}}$ at 
$(0, \ldots ,0,y)$ is full, 

with $$h _k (x_1, \ldots , x_{n+p(s)+2})= Q_k (s)
\big( g^k_1 (\overline{x}^{\sigma(s)^k_1}), \ldots , g^k_{q(s)} (\overline{x}^{\sigma (s) ^k_{q(s)}})\big),$$
and $$\overline{x}^{\sigma (s)^k_j} =(x_{\sigma(s)^k_j (1)} , \ldots ,x_{\sigma (s)^k_j (n)}).$$
\end{enumerate}

Then by Implicit Function Theorem, we have a mapping  
$$\Phi ^s :  M_s \longrightarrow \R \{ Y_1, \ldots , Y_{n+1} \}$$
which, sends 
$$\bigg(y,(g^k _j)_{\substack{1\leq j \leq q(s) \\ 1\leq k \leq p(s) +1}} \bigg)$$
to the analytic function $f$ defined in a neighbourhood of $\overline{0}\in \R ^{n+1}$, satisfying 
\begin{itemize}
\item $f(0, \ldots , 0)=y$
\item there are analytic functions $(f_1, \ldots ,f_{p(s)})$ in a neighbourhood of $(0,\ldots ,0)$ such that the graph of 
$(f_1, \ldots , f_{p(s)},f)$ is, in a neighbourhood of $(0, \ldots , 0,y)$, the zero-set of the $h_i$'s.
\end{itemize}

\begin{rem}    
By the definition of $n$-regularity, if $f:U \to \R$ is $n$-regular at $\overline{0}\in \R ^{n+1}$ then the germ of 
$f$ at $\overline{0}$ is in $\bigcup _{s \in \N} \Phi ^s (M_s)$.    
\end{rem}

Let's denote by $\R_{D,E} [\, X_1, \ldots X_m]$ the set of polynomials  
in $k$ variables, of degree $<D$ and of order $\geq d$ at the origin, with real coefficients. 

\begin{defi}
We denote the truncation mapping by $$ \left. \begin{array}{rccl}
 T^m _{DE}: & \R \{ X_1, \ldots , X_{m} \} &\to & \R_{D,E} [\, X_1, \ldots X_m] \\
								  \\ 
	&h &\mapsto &\displaystyle{\sum _{D \leq |\nu | < E}  \frac{\partial ^{|\nu|}h}
{\partial \overline{X}^{\, \nu}}(\overline{0})\cdot \overline{X}^{\, \nu}} \end{array} \right. . $$
\end{defi}

The chain derivation rule and an easy induction on $E$ gives us the next proposition, which will be useful to 
deduce non-surjectivity of the $\Phi ^s$ from the non-surjectivity of some rational mapping $\Phi _{DE} ^s$ 
between finite dimensional spaces.

\begin{prop}
\label{DiagComm}
Given three integers $s$, $D$ and $E$ with $D<E$, let $\widetilde{M}_s$ be the image of $M_s$ 
by the truncation $\Pi :={\mathrm Id}\, \otimes {({T^n _{0E}} ^{\otimes q(s)})}^{\otimes (p(s)+1)}$ 
of power series: 
$$\Pi :\R \times {(\R \{ X_1, \ldots , X_n \} ^{q(s)})}^{p(s)+1} \to \R \times {(\R_{0,E} [\, X_1, \ldots , X_n ] ^{q(s)})}^{p(s)+1}.$$
 
Then there is a rational mapping $\Phi _{DE} ^s$ such that the following diagram

\centerline{
\xymatrix@L=2pt{
M_s			\ar[r]^-{\Phi ^ s} \ar[d]_{\Pi}	
				&\R \{ Y_1,\ldots ,Y_{n+1} \} \ar[d]^{T^{n+1} _{DE}}\\
\widetilde{M}_s  	\ar[r]_-{\Phi _{DE} ^s}
				& \R_{D,E} [\, Y_1, \ldots , Y_{n+1}] }}
is commutative.
\end{prop}

This proposition simply says that the derivatives at the origin of order $< E$ of an element $\xi$ in the image of 
$\Phi ^s$ depend only on $y$ and on the derivatives at the origin of order $< E$ of $g_j ^k$ 
used to define $\xi$ in the source space of $\Phi ^s$, and this in a rational manner.

\section{Translation in the source space}

\label{lostintranslation}

The previous section would help us produce, by a diagonal argument, an analytic function which is outside of the image of each 
$\Phi ^s$ and thus is not $n$-regular at $\overline{0}\in \R ^{n+1}$.

But what we want is a function which is \emph{nowhere} $n$-regular in a neighbourhood of $\overline{0}$.
Hence we have to look at $\overline{x} \mapsto h(\overline{\alpha} + \overline{x})$ as 
$\overline{\alpha}$ ranges in a neighbourhood (let's say $]-1,1[^{n+1}$) of $\overline{0}$; 
unfortunately, we lose the finite dimensional dependency we found in the previous section.
 
More precisely, for $\overline{\alpha} \in ]-1,1[^{n+1}$, if we let $\tau_{\overline{\alpha}}$ be the function that 
assigns to an analytic function $h$ near $[-1,1]^{n+1}$ the function $\overline{x} \mapsto h(\overline{x} + \overline{\alpha})$
(which is analytic near $\overline{0}$) 
, we don't have the equality 
$$T_{DE} ^{n+1} (\tau_{\overline{\alpha}} (h)) = 
T_{DE} ^{n+1} \big( \tau_{\overline{\alpha}}\, \big(T_{DE} ^{n+1} (h)\big) \big); $$
each partial derivative of $h_{\overline{\alpha}}$ at the origin depends on \emph{all} 
partial derivatives of $h$ at zero.

The aim of this section is to show that this dependency can however be handled by metric arguments.

We first equip each $\R_{D,E} [\, Y_1, \ldots , Y_{n+1} ]$ with the norm: 
$$\| \sum _{\nu} a_{\nu} \overline{Y} ^{\nu} \| _{\infty}= \max _{\nu} \{ |a_{\nu}| \} .$$

\begin{prop}
Let $(D_k)$ be a increasing sequence of integers, $\eta$ a positive real number, 
$\overline{\alpha}$ a point in $]-1,1[^{n+1}$, $h$ an analytic function in a neighbourhood of $[-1,1]^{n+1}$
and $K$ an integer.

If for all $k>K$, we have $$\| T^{n+1}_{D_k D_{k+1}} (h) \|_{\infty} \leq \frac{\eta}{2^k {(D_{k+1}!)}^{\, n+1}},$$ then 
$$\|T _{D_KD_{K+1}} ^{n+1} (\tau_{\overline{\alpha}}\, (h)) - 
T _{D_KD_{K+1}} ^{n+1} \big( \tau_{\overline{\alpha}}\, 
\big(T _{D_KD_{K+1}} ^{n+1} (h)\big) \big)\|_{\infty}\leq \eta. $$ 

\end{prop}
\label{petitreste}
This is an easy consequence of the fact that, if $D_k \leq |\mu |< D_{k+1}$, then 
$$\frac{\partial ^{|\mu |} \big(\tau _{\alpha}( h)\big)}{\partial Y_1 ^{\mu _1} \ldots Y_{n+1} ^{\mu _{n+1}}}
(\overline{0}) =
\sum _{j\geq k} \sum _{\substack{\nu_i\geq \mu _i\\ D_j\leq |\nu |<D_{j+1}}}
\frac{\partial ^{|\nu |}h}{\partial Y_1 ^{\nu _1} \ldots Y_{n+1} ^{\nu _{n+1}}}(\overline{0}) 
\cdot \prod _i \binom{\nu _i}{\mu _i}{\alpha _i} ^{\nu _i - \mu _i}
 $$
and $\displaystyle |\prod _i \binom{\nu _i}{\mu _i}  {\alpha _i} ^{\nu _i - \mu _i}| \leq {(D_{k+1} !)}^{n+1}$ if 
$|\nu| < D_{k+1}$ and $| \overline{\alpha} |\leq 1$.

\begin{rem}
\label{S_r}
The linear mapping 
$L _{\overline{\alpha}} ^k$ on $\R _{D_k , D_{k+1}} [ Y_1, \ldots , Y_{n+1} ]$
defined by 
$L _{\overline{\alpha}} ^k (P)= T_{D_k D_{k+1}} (\tau _{\overline{\alpha}} \, (P))$ is an isomorphism, since 
the image of a monomial $\overline{X}^{\nu}$ is the sum of $\overline{X}^{\nu}$ and some lower degree monomials.  

Furthermore we have the identity 
$$\|{(L_{\overline{\alpha}} ^k)}^{-1} \| _ {\infty} = 
\max \{ 1/\| L _{\overline{\alpha}} ^k (P) \|_{\infty} \, ; \| P \| _ {\infty} =1 \}$$ and the mapping 
$(P, \overline{\alpha} ) \mapsto  1/\| L _{\overline{\alpha}} ^k (P)\| _{\infty}$ is continuous on the compact set  
$\{\| P \| _ {\infty} =1 \} \times {[-1,1]}^{n+1}$.

Thus we have a bound $S_k$ for the norm of 
${(L_{\overline{\alpha}} ^k)}^{-1}$, independent on $\overline{\alpha}\in ]-1,1[^{n+1}$.
\end{rem}


\section{Construction}
\label{construction}

We will use the good behaviour through truncation of the $\Phi ^s$ to build a sequence of integers $(D_s)$ and, for each 
$s \in \N$, a polynomial $P_s$ in $\R_{D_s , D_{s+1}}[ \, Y_1, \ldots Y_n ]$,
such that the formal power $h(Y_1, \ldots, Y_{n+1})=\sum _ s P_s (Y_1, \ldots , Y_{n+1})$ is the power expansion of an analytic 
function on $\R^{n+1}$ but such that $\tau _{\overline{\alpha}}\, (h)$ is outside of the image of $\Phi ^ s$ for each $s\in \N$ and 
$\overline{\alpha} \in {]-1,1[} ^{n+1}$. 
The restriction to ${[-1,1]}^{n+1}$ of this function (which is clearly $\R_{\mathrm{an}(n+1)}$-definable), 
will thus not be $\R_{\mathrm{an}(n)}$-definable as announced in section \ref{enumeration}.  

As we noted before Proposition \ref{DiagComm}, the lack of surjectivity of each $\Phi^s$ will follow from the lack of surjectivity of 
some mapping $\Phi _{D_s D_{s+1}} ^s$ between finite dimensional spaces. 

More precisely, if we fix a $s$ and a $D$, the function 
$$E \mapsto \dim (\R \times { 
(\R_{0,E} [\, X_1, \ldots , X_n ] ^{q(s)})}^{p(s)+1})$$ is a polynomial of degree $n$ in $E$, whereas
$$E \mapsto \dim (\R_{D,E} [\, Y_1, \ldots , Y_{n+1}])$$ is a polynomial of degree $n+1$.

We thus can build an increasing sequence of integers $(D_s)$ such that 
$$\dim (\R \times {( 
\R_{0,D_{s+1}} [\, X_1, \ldots , X_n ] ^{q(s)})}^{p(s)+1}) + n+1 $$ is smaller than 
$$\dim (\R_{D_s , D_{s+1}} [ \, Y_1, \ldots , Y_{n+1}]),$$ for each $s$.

Suppose we have built for $r<s$ some $P_r \in \R_{D_r,D_{r+1}}[\, Y_1, \ldots Y_{n+1} ]$ and $\eta _r >0$ such 
that \begin{enumerate}
\item[($A_r$):] $\forall t<r \, , \, \displaystyle \| P_r \| _{\infty} \leq \frac{\eta _t}{2^r {({D_{r+1}!})}^{n+1} }$ 
\item[($B_r$):] the ball of center $P_r$ and radius $\eta _r S_r$ 
(where $S_r$ is such that $\forall \overline{\alpha} \in ]-1,1[ ^{n+1},\, S_r \geq \| {(T _{D_rD_{r+1}} ^{n+1} \circ \tau _{\overline{\alpha }})}^{-1} \|_{\infty} $; 
see remark \ref{S_r})
does not meet the image of 
$\rho _r: (\alpha ,\xi) \mapsto {(T _{D_rD_{r+1}} ^{n+1} \circ \tau_{\overline{\alpha}})}^{-1} \circ \Phi _{D_r D_{r+1}} ^r (\xi)$ 
(where $\alpha $ ranges in $]-1,1[^{n+1}$ and $\xi$ in $\widetilde{M} _s$).
\end{enumerate}

We can then chose a $P_s\in \R_{D_s,D_{s+1}} [\, Y_1, \ldots Y_{n+1} ] $ and $\eta_s >0$ 
satisfying ($A_s$) and ($B_s$):

let $ \displaystyle \delta = \min \{\frac{\eta _t}{2^r {({D_{r+1}!})}^{n+1}}\, ; t<s \}$; by dimensional inequality  of 
source and image space (due to the choice of $D_{s+1}$) and rationality of 
$\rho _s : (\alpha ,\xi) \mapsto {(T _{D_sD_{s+1}} ^{n+1} \circ \tau_{\overline{\alpha}})}^{-1} \circ \Phi ^s _{D_s D_{s+1}} (\xi)$, 
we know that the image of $\rho _s$ is nowhere dense in $ \R _{D_s,D_{s+1}} [Y_1, \ldots Y_{n+1} ] $. 

We thus can 
chose a $P_s$ and $\eta _s$ such that $\| P_s \|< \delta$ and 
$$B(P_s \, ; \eta _s S_s) \cap \rho _s ({]-1,1[}^{n+1} \times \widetilde{M}_s) = \emptyset.$$

Let $h(Y_1, \ldots , Y_{n+1} ) $ be the formal series $\sum _s P_s (Y_1, \ldots , Y_{n+1})$.

We easily get from the inequalities ($A_r$) that $h$ is the power expansion of an analytic function on $\R ^{n+1}$.

Let $\overline{\alpha }$ be a point in ${]-1,1[}^{n+1}$.

From condition ($B_r$) we get that 
$$(T _{D_rD_{r+1}} ^{n+1} \circ \tau_{\overline{\alpha}}) 
\big( B(T _{D_rD_{r+1}} ^{n+1} h \, ; \, \eta _r S_r) \big)
\cap 
T _{D_rD_{r+1}} ^{n+1} \Phi _r (M_r) = \emptyset $$ 
and then by definition of $S_r$,  
$$  B((T _{D_rD_{r+1}} ^{n+1} \circ \tau_{\overline{\alpha}} \circ T _{D_rD_{r+1}} ^{n+1})\, h \, ; \, \eta _r)
\cap 
T _{D_rD_{r+1}} ^{n+1} \Phi _r (M_r) = \emptyset. $$

By ($A_s$) for $s>r$, we get from Proposition \ref{petitreste} that  
$$\| (T _{D_rD_{r+1}} ^{n+1} \circ \tau_{\overline{\alpha}})\, h - 
(T _{D_rD_{r+1}} ^{n+1} \circ \tau_{\overline{\alpha}} \circ T _{D_rD_{r+1}} ^{n+1})\, h \|_{\infty} \leq \eta _r;$$ 
thus $$  (T _{D_rD_{r+1}} ^{n+1} \circ \tau_{\overline{\alpha}})\,  h \notin T _{D_rD_{r+1}} ^{n+1} \Phi _r (M_r).$$

Hence $$ \tau_{\overline{\alpha}}\,  h \notin  \Phi _r (M_r).$$

\section*{Acknowledgement}

The author wishes to express his gratitude to { D. J.  Miller} 
(the { Fields Institute at Toronto}) for proving Proposition \ref{Miller} 
and to { R. Soufflet} (Universit\'e Lyon I) for suggesting the problem.

Many thanks also to { K.  Kurdyka},
{ P.  Speissegger} (and the logic team of McMaster University), 
and { G.  Valette} for all the suggestions. 


\end{document}